\documentclass[12pt]{article}

\DeclareMathAlphabet\mathbfcal{OMS}{cmsy}{b}{n}

\usepackage{geometry}
\usepackage[dvipsnames]{xcolor}
\usepackage{enumerate}
\usepackage[utf8]{inputenc}
\usepackage{amsfonts}
\usepackage{amsthm}
\usepackage{lscape}
\usepackage{mathrsfs}
\usepackage{graphicx}
\usepackage{tikz-cd}
\usepackage{calc}
\usepackage{float}
\usepackage{multirow}
\usepackage{lipsum}
\usepackage{amssymb}
\usepackage{mathtools}
\usepackage{mathdots}
\usepackage{fancyhdr}
\usepackage{tikz}
\usepackage{tikz-cd}
\usepackage{faktor}
\usepackage{setspace}
\usetikzlibrary{decorations.markings}
\newtheorem{theorem}{Theorem}
\theoremstyle{plain}
\usepackage{sectsty}
\usepackage{kpfonts}
\usepackage{blindtext}
\usepackage[all]{xy}
\usepackage{scalerel}

\newtheorem{lemma}{Lemma}

\newtheorem*{observation*}{Observation}
\newtheorem*{proposition*}{Proposition}
\newtheorem*{theorem*}{Theorem}
\newtheorem*{claim*}{Claim}  

\newcommand*{\QEDB}{\null\nobreak\hfill\ensuremath{\Box}}
\newcommand{\End}{{\mbox{End}}}
\newcommand{\GL}{{\mbox{GL}}}
\newcommand{\Aut}{\mbox{Aut}}

\title{A technical remark on the Donaldson-Futaki invariant for Fano reductive group compactifications}
\date{}
\author{Gabriella Clemente}
\begin{document}
\maketitle

\begin{abstract}
We present an elementary way of recovering a well-known criterion of K-stability for Fano reductive group compactifications.
\end{abstract}
 
\section{Introduction}\label{1}
A metric on a manifold is said to be canonical if its curvature is optimal in some variational sense. In K{\"a}hler geometry, the best studied canonical metrics are the extremal ones, which are critical points of the Calabi functional. Extremal metrics include the constant scalar curvature K{\"a}hler (cscK) metrics, among which are the K{\"a}hler-Einstein (KE) metrics. 

The Yau-Tian-Donaldson (YTD) conjecture is an umbrella term for the problem of finding necessary and sufficient algebro-geometric conditions for the existence of an extremal K{\"a}hler metric. See Theorem \ref{KE} for the YTD conjecture for positive KE metrics. 

In this note, we compute the Donaldson-Futaki invariant of a test-configuration of a polarized reductive group compactification using elementary methods (c.f.\ Theorem \ref{Futa}). Then we show how to apply that computation to recover a well-known combinatorial criterion of K-stability (Theorem \ref{tibi}). Let us begin with a brief review of the K{\"a}hler geometry that will be needed in the next section. For a thorough treatment of many of the topics that we present in sections \ref{1.1} and \ref{1.2} below, see \cite{Paul, Gabor}.

\subsection{The Mabuchi functional}\label{1.1}
Let $X$ be a compact K{\"a}hler manifold of dimension $n.$ Fix a K{\"a}hler class $\Omega \in H^{1,1}(X) \cap H^2(X,\mathbb{R}),$ and a basepoint $\omega_0 \in \Omega.$ Let \[\mathcal{M}_{\Omega}=\{\phi \in C^{\infty}(X,\mathbb{R}) \mid \omega_{\phi} :=\omega_0 + i \partial \bar{\partial} \phi >0\}\] be the space of K{\"a}hler potentials of K{\"a}hler metrics in $\Omega,$ $\overline{S}:=\frac{2n \pi c_1 (X) \cup \Omega^{n-1}}{\Omega^n}$ be the average scalar curvature, and $S_{\phi}$ be the scalar curvature of $\omega_{\phi}.$ The scalar curvature $1$-form $\sigma$ on $\mathcal{M}_{\Omega}$ is given as \[\sigma_{\phi}(\psi):=\langle \psi,\overline{S}-S_{\phi}\rangle_{L^2(\omega_{\phi})}=\int_X \psi (\overline{S}-S_{\phi}) \omega^n_{\phi},\] for any $\phi \in \mathcal{M}_{\Omega},$ $\psi \in T_{\mathcal{M}_{\Omega},\phi}=C^{\infty}(X,\mathbb{R}).$ 

The Mabuchi functional (or K-energy) is the uniquely defined functional $\mathcal{K}:\mathcal{M}_{\Omega} \to \mathbb{R}_{\geq 0},$ satisfying $d\mathcal{K}=\sigma$ and $\mathcal{K}(0)=0.$ In fact, $\mathcal{K}$ descends to a functional $\widehat{\mathcal{K}}$ on $\Omega$ whose Euler-Lagrange equation is the cscK equation, $S(\omega)=\overline{S},$ $\omega \in \Omega$ \cite{Mabuchi}. The K-energy has some interesting properties. For example, it is geodesically convex, and geodesics here are defined relative to the Riemannian metric on $\mathcal{M}_{\Omega},$ which at any $\phi \in \mathcal{M}_{\Omega}$ is given by $\langle f_1,f_2\rangle_{\phi}=\int_X f_1 f_2 \omega^n_{\phi},$ for any $f_1, f_2 \in T_{\mathcal{M}_{\Omega},\phi}.$ CscK metrics can be seen to be energy minimizers, although there may exist other minimizers as well.

\subsection{K-stability}\label{1.2}
Let $L \to X$ be an ample line bundle. Then for some $r>0,$ there is an embedding $\epsilon:X \hookrightarrow \mathbb{C}\mathbb{P}^{N_r},$ $\epsilon(p)=[s_0(p):\dots:s_{N_r}(p)],$ where $(s_i)_{i=0}^{N_r}$ is a basis of the space of global holomorphic sections $H^0(X,L^r).$ A test-configuration of the polarization $(X,L)$ is a 1-parameter deformation $(X_t,L_t)$ such that $X_t \simeq X$ unless possibly when $t=0.$ Precisely, it is the data of (1) a choice of embedding $\epsilon:X \hookrightarrow \mathbb{C}\mathbb{P}^{N_r}$ and (2) a $\mathbb{C}^*$-action on $\mathbb{C}\mathbb{P}^{N_r},$ or equivalently, a 1-parameter subgroup $\lambda:\mathbb{C}^* \to \GL_{N_r + 1}(\mathbb{C}).$

Now, let $\omega_{FS}$ be the Fubini-Study metric on $\mathbb{C}\mathbb{P}^{N_r}.$ For all $t \neq 0,$ $\omega_t:=\frac{1}{t} \epsilon^* \omega_{FS}$ is a K{\"a}hler metric on $X_t.$ Up to normalization, the Donaldson-Futaki (DF) invariant of the test configuration $(\epsilon,\lambda)$ is the rational number $DF(\epsilon,\lambda)=\lim_{t\to\infty} \frac{d\widehat{\mathcal{K}}(\omega_t)}{dt}.$ So the DF-invariant describes the asymptotics of Mabuchi's functional along special metric degenerations. The polarization $(X,L)$ is K-stable if the DF-invariant is positive for all test-configurations $(X_t,L_t)$ with $X_0 \not\simeq X.$ K-stability is a necessary and sufficient condition for the existence of a KE metric on a Fano manifold.

\begin{theorem}{(\cite{BM, CDS, Tian})}\label{KE}
Let $X$ be a compact K{\"a}hler Fano manifold. Then, there exists a KE form $\omega \in 2\pi c_1(X)$ iff $(X, K^{-1}_X)$ is K-stable, where $K^{-1}_X$ is the anti-canonical line bundle of $X.$ Moreover, whenever such a form exists, it is unique up to the standard action of the identity component $\Aut_0(X)$ of the group of holomorphic automorphisms of $X.$
\end{theorem}

The articles \cite{CDS, Tian} treat the existence part of the theorem, while \cite{BM} addresses the uniqueness part.

In practice, the main issue one encounters when trying to check the K-stability condition is classifying all inequivalent test-configurations. This issue can be overcome in highly symmetric cases, such as Fano reductive group compactifications, among which are Fano toric varieties.

\section{The Donaldson-Futaki invariant and reductive group compactifications}\label{2}

\subsection{Reductive group compactifications}\label{rgcs}
Good references for reductive groups and their representation theory are \cite{Borel1, French, reptheory, Sp}.

Let $G$ be a connected complex reductive group. A reductive group compactification (for $G$) is a normal irreducible projective $G \times G$-variety $X$ that has an open dense orbit that is equivariantly isomorphic to $G,$ where the $G \times G$-action on $G \subseteq X$ is $(g_1,g_2)x = g_1 x g_2^{-1}.$ In this context, a polarization is an ample $G \times G-$linearized line bundle on $X.$ Reductive group compactifications are a special instance of the stable reductive varieties from \cite{Brion}.  

Fix a choice of maximal torus $T \subset G$ with character lattice $M,$ and Lie algebra $\mathfrak{t}.$ Then, $\mathfrak{t}$ can be identified with $M_{\mathbb{R}}:=M \otimes \mathbb{R},$ and its dual $\mathfrak{t}^*$ with $N_{\mathbb{R}}:=N \otimes \mathbb{R},$ where $N$ is the lattice of 1-parameter subgroups of $T.$ Let $W$ be the Weyl group of $(G,T),$ and let $\Phi$ denote the root system of $(G,T)$ with selected positive roots $\Phi^+.$ We declare $2 \rho $ to be the sum of these positive roots. Relative to $\Phi^+ \subset N_{\mathbb{R}},$ the positive Weyl chamber is \[M_{\mathbb{R}}^{+} :=\{x \in M_{\mathbb{R}} \mid \langle \alpha, x \rangle \geq 0 \mbox{ for all }\alpha \in \Phi^+\},\] where $\langle \cdot, \cdot \rangle$ is the natural dual pairing. 

There is a one-to-one correspondence between lattice points $\lambda \in M_{\mathbb{R}}^+$ and irreducible $G-$representations $E_{\lambda},$ and to a lattice point $\lambda \in M_{\mathbb{R}}^+$ corresponds a $G \times G-$representation $End(E_{\lambda}).$ The dimension of $\End(E_{\lambda})$ is a polynomial 

\[\dim(\End(E_{\lambda}))=(\dim(E_{\lambda}))^2=H_d (\lambda) + H_{d-1} (\lambda) + \dots\] in $\lambda,$ and here $H_d$ stands for the degree $d$ homogeneous part of the polynomial $\dim(\End(E_{\lambda})),$ $H_{d-1}$ stands for the degree $d-1$ part, and so on. 

The combinatorial counterpart of a polarized $G$-compactification $(X,L)$ is a $W$-invariant convex polytope $P.$ Set $P^+:=P \cap M^+_{\mathbb{R}}.$ 

\subsection{The Fano condition}\label{tfc}
This short section is based on the article \cite{Ru}. At the polytope level, Fano is the condition that the distance between $2 \rho$ and any codimension one face of $P^{+}$ that does not meet the boundary of the positive Weyl chamber is equal to one. Denote the Zariski closure of $T$ in $X$ by $Z,$ which is a toric subvariety. When $X$ is Fano, the support function $v$ of $P$ is of the form $v=v_{K_{\mathbb{C}}}+v_Z,$ where $v_{K_{\mathbb{C}}}(x)=\langle 2 \rho , x \rangle $ for all $x$ in the positive Weyl chamber, $v_{K_{\mathbb{C}}} (wx)=v_{K_{\mathbb{C}}} (x)$ for all $w \in W,$ and $v_Z(x)=-g_{-K_Z}(-x),$ where $-g_{-K_Z}$ is the support function of the anti-canonical line bundle of the subvariety $Z \subset X.$ Since $P$ is also the polytope of $Z,$ the associated fan $\Sigma_P$ determines the toric subvariety $Z.$ From the theory of toric varieties, $-K_Z=\sum_{\rho \in \Sigma(1)} D_{\rho},$ where $\Sigma(1)$ is the set of $1$-dimensional cones of $\Sigma_P$ and $D_{\rho}$ is a prime torus invariant divisor on $Z.$ The support function $g_{-K_Z}$ has the property that $g_{-K_Z}(u_{\rho})=-1$ for all $\rho \in \Sigma(1),$ where $u_{\rho}$ is the minimal generator of the ray $\rho.$ In particular, if $a_i$ is the inward pointing normal to the $i$-th codimension one face of $P,$ $g_{-K_Z}(a_i)=-g_{-K_Z}(-a_i)=-1.$ Then, $v(a_i)=\langle a_i, 2\rho \rangle -1,$ and so the facet presentation of the polytope is
\[P=\{x \in M_{\mathbb{R}} | \langle a_i, x \rangle \geq \langle a_i, 2 \rho \rangle -1 \}.\] Consequently, the equation that defines the $i$-th boundary face of $P$ is $f_i(x)=\langle a_i, x-2 \rho \rangle + 1$ so that $f_i(2 \rho)=1.$ 

\subsection{Donaldson-Futaki invariant}\label{dfi}
Let $(X,L)$ be a polarized reductive group compactification, and $P$ be its polytope realization. A test-configuration of $(X,L)$ is given by a piecewise linear (PL), $W$-invariant, convex function $f$ on $P$ that is affine linear on $P^{+}.$ 

\begin{theorem}\label{maincalcalex}{(Theorem 3.3 \cite{Alex})} Let $f$ be a convex rational $W-$invariant PL function on $P.$ Then the DF-invariant of the corresponding test-configuration is given by the formula 

\[-F_1(f)=\frac{1}{2 \int_{P^{+}} H_d d \mu} \Big( \int_{\partial P^{+}} fH_{d} d \sigma +2 \int_{P^{+}} fH_{d-1} d\mu-a\int_{P^{+}} fH_d d\mu \Big),\] where \[a=\frac{\int_{\partial P^{+}} H_d d\sigma + 2\int_{P^{+}} H_{d-1} d\mu}{\int_{P^{+}} H_d d \mu}.\] 
\end{theorem} 

Here $d \mu$ is the Lebesgue measure restricted to $P,$ and the boundary measure $d \sigma$ is a positive measure on $\partial P$ that is normalized so that on each codimension one face, which is defined by an equation $l(x):=\langle a,x \rangle = c,$ $d \sigma \wedge dl = \pm d \mu$ holds.

In the sequel, we obtain a number of identities that together with the Fano condition will allow us to simplify the DF-invariant from Theorem \ref{maincalcalex}. 

\begin{theorem}\label{Futa}
Assume that $P$ satisfies the Fano condition. Let $Vol_{DH}(P^{+}):=\int_{P^{+}} H_d d \mu$ and $bar_{DH}(P^{+}):=\frac{1}{Vol_{DH}(P^{+})} \int_{P^{+}} x H_d d \mu$ be the volume, respectively the barycenter of $P^{+}$ w.r.t.\ the Duistermaat-Heckman (DH) measure. Let $2\rho$ be the sum of the selected positive roots. Then, the DF-invariant of the test-configuration represented by $f$ is
\begin{equation*}
\begin{split}
-F_1(f)&=\frac{1}{2 Vol_{DH}(P^{+})} \int_{P^{+}} \langle \nabla f, x-2\rho\rangle H_d d\mu \\
&= \frac{1}{2} \langle bar_{DH}(P^{+}) - 2 \rho, \nabla f \rangle,
\end{split}
\end{equation*}
\end{theorem}

Theorem \ref{Futa} can be interpreted as a reductive analog of Theorem C in \cite{spherical}. In essence, our proof of Theorem \ref{Futa} and of Delcroix's barycentric critetion, implements an adjustment of a technique first used in the toric setting \cite{ZZ, ZZ1} to transform boundary terms into interior terms via the divergence theorem. See \cite{LLZ} for an application of this technique to group compactifications. The next section explains how to use Theorem \ref{Futa} to recover Delcroix's barycentric criterion for K-stability (cf.\ Theorem \ref{tibi}). 

Choose once and for all an isomorphism $M_{\mathbb{R}} \simeq \mathbb{R}^n$ so that $P$ can be viewed as though contained in $\mathbb{R}^n.$ 

\begin{lemma}\label{identities}
Let $\Phi^{+}=\{\alpha_1,\dots,\alpha_r\},$ $c=\prod_{i=1}^r \langle \alpha_i,\rho \rangle^2,$ where $\rho=\frac{1}{2} \sum_{i=1}^r \alpha_i,$ and let $\{e_j\}_{j=1}^n$ be the standard basis of $\mathbb{R}^n.$ Then,
\begin{enumerate}
\item \[H_d(x)=\frac{1}{c} \prod_{i=1}^r \langle \alpha_i, x \rangle^2,\] 

\item \[H_{d-1}(x)=\frac{1}{c} \sum_{j=1}^r2 \langle \alpha_j,x \rangle \langle \alpha_j, \rho \rangle \langle \alpha_1, x \rangle^2 \dots \widehat{\langle \alpha_j , x \rangle}^2 \dots \langle \alpha_r, x \rangle^2,\]

\item \[\nabla H_d(x)=\frac{1}{c} \sum_{j=1}^n \big(\sum_{i=1}^n 2 \langle \alpha_i, x \rangle \langle \alpha_i, e_j \rangle \langle \alpha_1, x \rangle^2 \dots \widehat{\langle \alpha_j , x \rangle}^2 \dots \langle \alpha_r, x \rangle^2 \big)e_j,\]

\item $\langle \nabla H_d(x), \rho \rangle = H_{d-1}(x),$

\item $\langle \nabla H_d(x), x \rangle=2rH_d(x),$ and 

\item for any smooth function $f: P \rightarrow \mathbb{R},$  \[div((x-2\rho)fH_d)=\langle \nabla f, x-2 \rho \rangle H_d+(2r+n)fH_d-2fH_{d-1}.\]
\end{enumerate}
\end{lemma}

\begin{proof}
Let $E_x$ be an irreducible representation with highest weight $x.$ To prove 1.\ and 2., we make use of the Weyl dimension formula \[\dim(E_x)=\frac{\prod_{i=1}^r \langle \alpha_i,x+\rho \rangle}{\prod_{i=1}^r \langle \alpha_i,\rho \rangle}.\] From the expression

\[\dim(E_x)^2=\frac{1}{c} \prod_{i=1}^r (\langle \alpha_i, x \rangle^2 + 2 \langle \alpha_i, x \rangle \langle \alpha_i, \rho \rangle + \langle \alpha_i, \rho \rangle^2),\] it follows that if $d$ is the highest degree homogeneous part of the polynomial $\dim(E_x)^2,$ then \[H_d(x)=\frac{1}{c} \prod_{i=1}^r \langle \alpha_i, x \rangle^2,\] and the $(d-1)-$degree homogeneous part of $\dim(E_x)^2$ is  \[H_{d-1}(x)=\frac{1}{c} \sum_{j=1}^r2 \langle \alpha_j,x \rangle \langle \alpha_j, \rho \rangle \langle \alpha_1, x \rangle^2 \dots \widehat{\langle \alpha_j , x \rangle}^2 \dots \langle \alpha_r, x \rangle^2. \]

For 3., note that $\frac{\partial}{\partial x_j} \langle \alpha_i , x \rangle = \langle \alpha_i, e_j \rangle$ so that

\[\frac{\partial}{\partial x_j} H_d(x) = \frac{1}{c} \sum_{i=1}^r 2 \langle \alpha_i,x \rangle \langle \alpha_i , e_j \rangle \langle \alpha_1, x \rangle^2 \dots \widehat{\langle \alpha_i , x \rangle}^2 \dots \langle \alpha_r, x \rangle^2,\] and hence 

\begin{equation*}
\nabla H_d(x) = \sum_{j=1}^r \frac{\partial}{\partial x_j} H_d(x) e_j \\
 = \frac{1}{c} \sum_{j=1}^n \big(\sum_{i=1}^n 2 \langle \alpha_i, x \rangle \langle \alpha_i, e_j \rangle \alpha_1, x \rangle^2 \dots \widehat{\langle \alpha_j , x \rangle}^2 \dots \langle \alpha_r, x \rangle^2 \big)e_j.
\end{equation*}

For 4., notice that 

\begin{equation*}
\begin{split}
\nabla H_d(x) & = \frac{1}{c} \sum_{j=1}^n \big(\sum_{i=1}^n 2 \langle \alpha_i, x \rangle \langle \alpha_i, e_j \rangle \langle \alpha_1, x \rangle^2 \dots \widehat{\langle \alpha_j , x \rangle}^2 \dots \langle \alpha_r, x \rangle^2 \big)e_j \\
& = \sum_{i=1}^r \big( \frac{1}{c} 2 \langle \alpha_i, x \rangle \langle \alpha_1, x \rangle^2 \dots \widehat{\langle \alpha_i , x \rangle}^2 \dots \langle \alpha_r, x \rangle^2 \big) \sum_{j=1}^n \langle \alpha_i, e_j \rangle e_j \\
& = \sum_{i=1}^r \big( \frac{1}{c} 2 \langle \alpha_i, x \rangle \langle \alpha_1, x \rangle^2 \dots \widehat{\langle \alpha_i , x \rangle}^2 \dots \langle \alpha_r, x \rangle^2 \big) \alpha_i
\end{split}
\end{equation*}

and then 

\[\langle \nabla H_d(x), \rho \rangle = \sum_{i=1}^r \frac{1}{c} 2 \langle \alpha_i, x \rangle \langle \alpha_i, \rho \rangle \langle \alpha_1, x \rangle^2 \dots \widehat{\langle \alpha_i , x \rangle}^2 \dots \langle \alpha_r, x \rangle^2 = H_{d-1}(x).\]

For 5., observe that since \[\nabla H_d(x) = \sum_{i=1}^r \big( \frac{1}{c} 2 \langle \alpha_i, x \rangle \langle \alpha_1, x \rangle^2 \dots \widehat{\langle \alpha_i , x \rangle}^2 \dots \langle \alpha_r, x \rangle^2 \big) \alpha_i,\] and since for each $i,$ 
\[\big\langle \frac{1}{c} 2 \langle \alpha_i, x \rangle \langle \alpha_1, x \rangle^2 \dots \widehat{\langle \alpha_i , x \rangle}^2 \dots \langle \alpha_r, x \rangle^2 \alpha_i, x \big\rangle =  2\big(\frac{1}{c} \langle \alpha_1, x \rangle^2 \dots \dots \langle \alpha_r, x \rangle^2\big),\] indeed we have that

\[\langle \nabla H_d(x), x \rangle = 2r\Big(\frac{1}{c} \prod_{i=1}^r \langle \alpha_i,x\rangle^2\Big) = 2r H_d(x).\]

The above identities now imply the last point. Namely, 
\begin{equation*}
\begin{split}
div((x-2\rho)fH_d)&=\langle \nabla (fH_d), x-2\rho \rangle + div(x-2\rho)fH_d\\
& = \langle \nabla f, x-2\rho \rangle H_d+ \langle \nabla H_d , x-2\rho \rangle f+ n fH_d \\
& = \langle \nabla f, x-2\rho \rangle H_d+ \langle \nabla H_d, x \rangle f-2\langle \nabla H_d,\rho\rangle f+ n fH_d\\
& = \langle \nabla f, x-2\rho \rangle H_d+(2r+n)fH_d-2fH_{d-1}. 
\end{split}
\end{equation*}
\end{proof}

\noindent
\emph{Proof of Theorem \ref{Futa}}. Suppose that $\partial P^{+}$ has $k$ codimension one faces $\partial P_i^{+}.$ Let $\{\partial P_i^{+} : i=1,\dots,m\}$ be the set of all codimension one faces of $P^{+}$ that do not intersect the boundary of the positive Weyl chamber. Suppose that $\partial P_i^{+}$ is defined by $\langle a_i,  x \rangle - c_i=0$ and set $f_i(x):=\langle a_i,  x \rangle - c_i.$ The (inward) unit normal vector field to $\partial P_i^{+}$ is $-\frac{\nabla f_i}{\| \nabla f_i \|}=-\frac{a_i}{\|a_i\|}.$ Since $P$ satisfies the Fano condition, for $x \in \partial P_i^{+},$ we have that

\begin{equation*}
\begin{split}
\langle (x-2 \rho )f H_d, -\frac{a_i}{\|a_i\|} \rangle &=-H_d f \Big( \Big\langle x-2\rho , \frac{a_i}{\| a_i \|} \Big\rangle \Big) \\
& = \frac{-H_d(x)f}{\| a_i\|} \Big( \langle x, a_i \rangle - \langle 2 \rho , a_i \rangle \Big) \\
& =\frac{H_d(x)f}{\| a_i\|} \Big(\langle 2 \rho , a_i \rangle-c_i\Big) \\
& = \frac{H_d(x)f}{ \| a_i \|}.
\end{split}
\end{equation*}

The divergence theorem implies that \[\int_{P^{+}} div( (x-2 \rho )f H_d) d \mu = \sum_{i=1}^m \int_{ \partial P_i^{+}} \frac{H_d f}{\|a_i\|} d \sigma_i + \sum_{i=m+1}^k \int_{ \partial P_i^{+}} \langle (x-2 \rho )f H_d, -\frac{a_i}{\|a_i\|} \rangle d \sigma_i,\] where $d \sigma_i$ is the standard Lebesgue measure on $\partial P$ with domain restricted to $\partial P_i.$ When $i=m+1,\dots,k,$ $\partial P_i^{+}$ is in the boundary of the positive Weyl chamber, and \[\int_{ \partial P_i^{+}} \langle (x-2 \rho )f H_d, \frac{a_i}{\|a_i\|} \rangle d \sigma_i = 0.\] Then
\[\int_{P^{+}} div( (x-2 \rho )f H_d) d \mu = \sum_{i=1}^m \int_{ \partial P_i^{+}}  \frac{H_d(x)f}{ \| a_i \|} d \sigma_i \] and the right hand side is the definition of \[\int_{ \partial P^{+}} fH_d  d \sigma.\] 

By 6.\ of Lemma \ref{identities}, taking $f=1,$ we obtain that \[div((x-2\rho)H_d)=(2r+n)H_d-2H_{d-1}.\] Then, by the divergence theorem, \[\int_{\partial P^{+}} H_d d \sigma =(2r+n) \int_{P^{+}} H_d d\mu - 2 \int_{P^{+}} H_{d-1} d \mu.\] Hence, \[a=\frac{\int_{\partial P^{+}} H_d d\sigma + 2\int_{P^{+}} H_{d-1} d\mu}{\int_{P^{+}} H_d d \mu}=2r+n.\] 

Upon substituting the above calculations into Alexeev's and Katzarkov's DF-invariant (cf.\ Theorem \ref{maincalcalex}), again using 6.\ of Lemma \ref{identities} to rewrite the first integral, we find that
\[-F_1(f)=\frac{1}{2 \int_{P^{+}} H_d d \mu} \int_{P^{+}} \langle \nabla f, x-2\rho\rangle H_d d\mu.\] 

Suppose that $f$ on $P^{+}$ is given as $f(x)=\sum_{j=1}^n b_jx_j + k.$ Put $b=(b_1,\dots, b_n),$ $x=(x_1, \dots, x_n)$ and $2 \rho=(2 \rho_1, \dots, 2 \rho_n),$ and let $e_j$ be the $j-$th standard basis vector of $\mathbb{R}^n.$ Then, $\langle \nabla f, x - 2 \rho \rangle = \sum_{j=1}^n b_j(x_j - 2 \rho_j)$ and it follows that 

\begin{equation*}
\begin{split}
-F_1(f) &= \frac{1}{2 Vol_{DH} (P^{+})} \int_{P^{+}} \langle \nabla f, x - 2 \rho \rangle H_d d \mu \\
& = \frac{1}{2 Vol_{DH}(P^{+})} \Big(\sum_{j=1}^n b_j \int_{P^{+}} x_j H_d d \mu - \sum_{j=1}^n b_j (2 \rho_j) Vol_{DH}(P^{+}) \Big)\\
& = \frac{1}{2} \Big(\langle bar_{DH}(P^{+}), \sum_{j=1}^n b_j e_j \rangle - \sum_{j=1}^n b_j (2 \rho_j) \Big) \\
&=\frac{1}{2} \big(\langle bar_{DH}(P^{+}), b \rangle - \langle b, 2 \rho \rangle \big) \\
& = \frac{1}{2} \langle bar_{DH}(P^{+}) -2\rho , \nabla f \rangle.
\end{split}
\end{equation*}
\QEDB

\subsection{K-stability as a combinatorial criterion}\label{Knew}

Next we show how to use Theorem \ref{KE} alongside Theorem \ref{Futa} to verify Delcroix's barycentric criterion (Theorem \ref{tibi}). 

From the works \cite{Dat} or \cite{Lu}, it follows that all relevant test-configurations are defined by affine linear functions on $P^{+} \subset \mathbb{R}^n.$ So when $(X_P, K_{X_P}^{-1})$ is $K$-stable, $-F_1(f) \geq 0$ for any linear function $f:P^{+} \rightarrow \mathbb{R}.$ In particular, for the coordinate function $x^i(x):=\langle x, e_i \rangle,$ where $e_i$ is the $i$-th standard basis vector of $\mathbb{R}^n,$ $-F_1(x^i) \geq 0.$ Then since \[-F_1(x^i)=\frac{1}{2Vol_{DH}(P^{+})} \int_{P^{+}} (x^i(x)-x^i(2 \rho))H_d d \mu,\] $-F_1(x^i) \geq 0$ implies $x^i(b) \geq x^i(2 \rho),$ where $b$ is the barycenter of $P^{+}$ w.r.t.\ the DH measure. 

Conversely, suppose that $x^i(b) \geq x^i(2 \rho)$ so that $-F_1(x^i) \geq 0.$  Take any affine linear function $f(x)=\sum_{i=1}^n a_ix^i + b$ on $P^{+}$ and note in order to extend $f$ in a $W-$invariant way so that it is a convex function on $P,$ we need that $a_i \geq 0$ for all $i.$ Then since \[\langle \nabla f, x-2\rho \rangle =\sum_{i=1}^n a_i(x^i(x)-x^i(2 \rho)),\]
\[-F_1(f)=\sum_{i=1}^n a_i(-F_1(x^i)) \geq 0,\] and so the anti-canonical polarization of $X_P$ is K-stable. 

\begin{theorem}{(Theorem 1.5 \cite{tibiphd}, Theorem A \cite{Del})}\label{tibi}
Let $X$ be a smooth Fano compactification of a reductive group $G$ with associated polytope $P.$ Let $\Phi$ be the root system of a pair $(G,T).$ The barycenter of $P^{+}$ w.r.t.\ the DH measure is $b:=\frac{1}{\int_{P^{+}} H_d d \mu} \int_{P^{+}} x H_d d \mu.$ There exists a KE metric on $X$ iff $b \in 2 \rho + \Theta,$ where $\Theta$ is the relative interior of the cone generated by $\Phi^+.$
\end{theorem}

Theorem \ref{tibi} has been generalized to smooth Fano spherical varieties (cf.\ Theorem A \cite{spherical}). 

Indeed, Theorem \ref{tibi} implies Wang's and Zhu's barycenter characterization of Fano KE toric varieties \cite{toribari} (see also \cite{mabi}). Namely, $X_P$ is KE iff the barycenter of $P$ is the origin.\\

\noindent
\textbf{Acknowledgement.} I thank Richard Hind and Jean-Pierre Demailly$^{\dag}$ for once encouraging me to make this work available.

\begin{minipage}[t]{10cm}
\begin{flushleft}
\small{
\textsc{G.\ Clemente}
\\*D{\'e}partement de Math{\'e}matiques, CNRS UMR 6134 SPE
\\*Universit{\`a} di Corsica Pasquale Paoli
\\*Corte, 20250, France 
\\*e-mail: clemente\_g@univ-corse.fr
}
\end{flushleft}
and
\begin{flushleft}
\small{
\textsc{}
\\*IHES, CNRS UMR 9009
\\*Bures-sur-Yvette, 91440, France
\\*e-mail: clemente@ihes.fr
}
\end{flushleft}
\end{minipage}

\end{document}